# A method for designing involute trajectory timelike ruled surfaces in Minkowski 3-space

**Mustafa Bilici**

Ondokuz Mayıs University, Education Faculty,

Department of Mathematics, 55200 Samsun – Turkey

**Abstract**

The aim of this paper is to present a new perspective on the generation of developable trajectory ruled surfaces in Minkowski 3-space. Involute trajectory ruled surfaces generated by the Frenet trihedron, moving along spacelike involutes of a given timelike space curve, is stated according to Lorentzian timelike angle between tangent vector and unit vector of direction Darboux vector of this timelike space curve. Also some new results and theorems related to the developability of the involute trajectory timelike ruled surfaces are obtained. Finally we illustrate these surfaces by presenting one example.

e-mail: mbilici@omu.edu.tr (M Bilici)

## 1. Introduction

In classical differential geometry, ruled surfaces and their special subclass, developable surfaces are well known concepts. Ruled surfaces can be described as the set of points swept by a moving straight line. For example, a cone is formed by keeping one point of a line fixed whilst moving another point along a circle.

A developable ruled surface is a surface that can be isometrically mapped (i.e. developed) into the plane. Among the ruled surfaces, the developables play a distinguished role since they have natural applications in many areas of engineering and manufacturing. For instance, an aircraft designer uses them to design the airplane wings, and a tinsmith uses them to connect two tubes of different shapes with planar segments of metal sheets.

In the past, offsets of ruled surfaces have been subject of many studies: In [13], Ravani and Ku generalized the theory of Bertrand curves for Bertrand ruled surface-offsets, using line geometry. Some of the studies in the Computer Aided Geometric Design (CAGD) literature dealing with offsets of surfaces have been given by [2, 3, 4 and 13]. The corresponding characterizations of timelike and spacelike ruled surfaces in Minkowski 3-space have been given by Turgut and Hacısalihoğlu [6-9]. Yaylı [11] obtained the distribution parameter of a spacelike ruled surface under the one-parameter spatial motion and proved that the spacelike ruled surface is developable if and only if the base curve is helix. Woestijne [14] and Young Ho Kim and Dae Won Yoon [15] classified the Lorentz surfaces. These studies are important for space kinematic and mechanism. Furthermore, the geometry of trajectory ruled surfaces is widely applied to the study of design problems in spatial mechanisms or space kinematics. In addition to all these references, for understanding the space kinematics and mechanism, references [1, 5, 10] and references quoted therein are advised for readers interested in this field.

In this paper, the generation of trajectory ruled surfaces in Minkowski 3-space generated by the Frenet trihedron moving along involutes of a given timelike space curve is stated according to Lorentzian timelike angle between the binormal vector and unit vector of direction Darboux vector of this space curve. And some new results and theorems related to the developability of involute trajectory ruled surfaces are obtained. Additionally these surfaces have a number of applications in computer-aided design.

**2. Preliminaries**

Let Minkowski 3-space $\mathbb{R}_1^3$ be the vector space $\mathbb{R}^3$ provide with the Lorentzian inner product,

$$\langle \vec{u}, \vec{v} \rangle = -u_1 v_1 + u_2 v_2 + u_3 v_3,$$

where $\vec{u} = (u_1, u_2, u_3), \vec{v} = (v_1, v_2, v_3) \in \mathbb{R}^3$.

A vector $\vec{u} \in \mathbb{R}^3$ is said to be timelike if $\langle \vec{u}, \vec{u} \rangle < 0, \vec{u} \neq 0$ spacelike if $\langle \vec{u}, \vec{u} \rangle > 0$ and lightlike (or null) if $\langle \vec{u}, \vec{u} \rangle = 0$ for all. Similarly, an arbitrary curve $\vec{r} : I \to \mathbb{R}_1^3,\ I \in \mathbb{R}$ in $\mathbb{R}_1^3$ can locally be timelike, spacelike or null (lightlike), if all of its velocity vectors $\dot{\vec{r}}(s)$ are timelike, spacelike or null, for every $s \in I \subset \mathbb{R}$, respectively [12]. If $\langle \vec{u}, \vec{v} \rangle = 0$ for all $\vec{u}$ and $\vec{v}$, the vectors $\vec{u}$ and $\vec{v}$ are called perpendicular in sense of Lorentz. The norm of $\vec{u} \in \mathbb{R}_1^3$ is defined as $\|\vec{u}\| = \sqrt{|\langle \vec{u}, \vec{u} \rangle|}$ .A timelike vector $\vec{u}$ is said to be positive (resp. negative)if and only if $u_1 > 0$ (resp. $u_1 < 0$). The vector product of $\vec{u}$ and $\vec{v}$ given by [16] as:

$$\vec{u} \wedge \vec{v} = (u_3 v_2 - u_2 v_3, u_1 v_3 - u_3 v_1, u_1 v_2 - u_2 v_1). \tag{2.1}$$

A surface in $\mathbb{R}_1^3$ is called timelike surface if the induced metric on the surface is Lorentz metric, i.e. the normal on the surface is a spacelike vector. Let $\vec{r}(s)$ be a unit speed timelike space curve with curvature $\kappa$ and torsion $\tau$. We denote by $\{\vec{t}(s), \vec{n}(s), \vec{b}(s)\}$ the moving Frenet frame along the curve $\vec{r}(s)$. Then $\vec{t}$, $\vec{n}$ and $\vec{b}$ are the tangent, the principal normal and the binormal vector of the curve $\vec{r}(s)$, respectively. In this trihedron, $\vec{t}$ is timelike vector, $\vec{n}$ and $\vec{b}$ are spacelike vectors. So scalar product and cross product of these vectors is given by

$$-\langle \vec{t},\vec{t}\rangle = \langle \vec{n},\vec{n}\rangle = \langle \vec{b},\vec{b}\rangle = 1, \quad \langle \vec{t},\vec{n}\rangle = \langle \vec{n},\vec{b}\rangle = \langle \vec{b},\vec{t}\rangle = 0, \tag{2.2}$$

$$\vec{t}\times\vec{n} = -\vec{b} \ , \ \vec{n}\times\vec{b} = \vec{t} \ , \ \vec{b}\times\vec{t} = -\vec{n}. \tag{2.3}$$

In this situation, the Frenet formulas are given by [14]

$$\dot{\vec{t}} = \kappa\vec{n}, \qquad \dot{\vec{n}} = \kappa\vec{t} - \tau\vec{b}, \qquad \dot{\vec{b}} = \tau\vec{n}. \tag{2.4}$$

We also define the rotation vector of this trihedron, the Darboux vector of $\vec{r}(s)$, which provides a useful way of interpreting $\kappa$ and $\tau$ geometrically. We can express this in terms of the Frenet-Serret vectors as [18]

$$\vec{D} = \tau\vec{t} - \kappa\vec{b} \tag{2.5}$$

There are two cases corresponding to the causal characteristic of Darboux vector $\vec{D}$.

**Case I.** If $|\kappa| > |\tau|$ then $\vec{D}$ is a spacelike vector. In this situation, we can write

$$\begin{cases} \kappa = \|\vec{D}\|\cosh\theta \\ \tau = \|\vec{D}\|\sinh\theta \end{cases}, \quad \|\vec{D}\|^2 = g(\vec{D},\vec{D}) = \kappa^2 - \tau^2 \tag{2.6}$$

and the unit vector $\vec{c}$ of direction $\vec{D}$ is

$$\vec{c} = \frac{1}{\|\vec{D}\|}\vec{D} = \sinh\theta\,\vec{t} - \cosh\theta\,\vec{b} \tag{2.7}$$

**Case II.** If $|\kappa| < |\tau|$ then $\vec{D}$ is a timelike vector. In this situation, we have

$$\begin{cases} \kappa = \left\| \vec{D} \right\| \sinh\theta \\ \tau = \left\| \vec{D} \right\| \cosh\theta \end{cases}, \quad \left\| \vec{D} \right\|^2 = -g\left(\vec{D}, \vec{D}\right) = \tau^2 - \kappa^2 \tag{2.8}$$

and the unit vector $\vec{c}$ of direction $\vec{D}$ is

$$\vec{c} = \frac{1}{\left\| \vec{D} \right\|} \vec{D} = \cosh\theta\, \vec{t} - \sinh\theta\, \vec{b}, \tag{2.9}$$

where $\theta$ is the Lorentzian timelike angle between $\vec{t}$ and $\vec{c}$.

**Remark 1.** We can easily see from equations of the section Case I and Case II that $\frac{\tau}{\kappa} = \tan\text{h}\,\theta \left( \text{or } \frac{\tau}{\kappa} = \coth\theta \right)$ and if $\theta$ is an arbitrary constant then $\vec{r}(s)$ is a general helix.

The angle between two vectors in Minkowski 3-space is defined by [19]:

**Definition 1.** Let $\vec{u}$ and $\vec{v}$ be spacelike vectors in $\mathbb{R}_1^3$ that span a spacelike vector subspace, then we have $\left| g\left(\vec{u}, \vec{v}\right) \right| \le \left\| \vec{u} \right\| \left\| \vec{v} \right\|$ and hence, there is a unique real number $\varphi$ such that

$$g\left(\vec{u}, \vec{v}\right) = \left\| \vec{u} \right\| \left\| \vec{v} \right\| \cos\varphi. \tag{2.10}$$

The reel number $\varphi$ is called the Lorentzian spacelike angle between $\vec{u}$ and $\vec{v}$.

**Definition 2**. Let $\vec{u}$ and $\vec{v}$ be spacelike vectors in $\mathbb{R}_1^3$ that span a timelike vector subspace, then we have $g\left(\vec{u}, \vec{v}\right) > \left\| \vec{u} \right\| \left\| \vec{v} \right\|$ and hence, there is a unique positive real number $\varphi$ such that

$$\left| g\left(\vec{u}, \vec{v}\right) \right| = \left\| \vec{u} \right\| \left\| \vec{v} \right\| \cosh\varphi. \tag{2.11}$$

The reel number $\varphi$ is called the Lorentzian timelike angle between $\vec{u}$ and $\vec{v}$.

**Definition 3.** Let $\vec{u}$ be a spacelike vector and $\vec{v}$ a positive timelike vector in $\mathbb{R}_1^3$, then there is a unique non-negative reel number $\varphi$ such that

$$\left|g\left(\vec{u},\vec{v}\right)\right| = \|\vec{u}\|\|\vec{v}\|\sinh\varphi . \tag{2.12}$$

The reel number $\varphi$ is called the Lorentzian timelike angle between $\vec{u}$ and $\vec{v}$.

**Definition 4.** Let $\vec{u}$ and $\vec{v}$ be positive (negative) timelike vectors in $\mathbb{R}_1^3$, then there is a unique non-negative real number $\varphi$ such that

$$g\left(\vec{u},\vec{v}\right) = \|\vec{u}\|\|\vec{v}\|\cosh\varphi . \tag{2.13}$$

The reel number $\varphi$ is called the Lorentzian timelike angle between $\vec{u}$ and $\vec{v}$.

Let $\vec{\gamma}(s)$ be spacelike involutes of a given timelike space curve $\vec{r}(s)$ in $\mathbb{R}_1^3$. If $\vec{\gamma}(s)$ denotes the involutes of $\vec{r}(s)$, we have

$$\vec{\gamma}(s) = \vec{r}(s) + (c-s)\vec{t}(s), \tag{2.14}$$

where $c$ is an arbitrary constant and $\dot{\vec{r}}(s) = \vec{t}(s)$. Let Frenet frames of $\vec{r}(s)$ and $\vec{\gamma}(s)$ be $\left\{\vec{t}(s),\vec{n}(s),\vec{b}(s)\right\}$ and $\left\{\vec{t}^*(s),\vec{n}^*(s),\vec{b}^*(s)\right\}$, respectively. Consider the Frenet frame $\left\{\vec{t}^*(s),\vec{n}^*(s),\vec{b}^*(s)\right\}$ attached to the spacelike involutes $\vec{\gamma}(s)$ such that the tangent vector $\vec{t}^*(s)$ is a unit spacelike vector, the principal normal vector $\vec{n}^*(s)$ is a unit timelike vector and the binormal vector $\vec{b}^*(s)$ is a unit spacelike vector. For that

$$\left\langle \vec{t}^*,\vec{t}^*\right\rangle = -\left\langle \vec{n}^*,\vec{n}^*\right\rangle = \left\langle \vec{b}^*,\vec{b}^*\right\rangle = 1, \quad \left\langle \vec{t}^*,\vec{n}^*\right\rangle = \left\langle \vec{n}^*,\vec{b}^*\right\rangle = \left\langle \vec{b}^*,\vec{t}^*\right\rangle = 0 \tag{2.15}$$

can be written. For the Frenet trihedron, the vectoral product is given by

$$\vec{t}^* \wedge \vec{n}^* = -\vec{b}^*, \ \vec{n}^* \wedge \vec{b}^* = -\vec{t}^*, \ \vec{b}^* \wedge \vec{t}^* = \vec{n}^* . \tag{2.16}$$

On the other hand, the relationships between the Frenet frames of $\vec{r}(s)$ and $\vec{\gamma}(s)$ can be given by [17]

1) If $\vec{D}$ is a spacelike vector, then

$$\begin{pmatrix} \vec{t}^*(s) \\ \vec{n}^*(s) \\ \vec{b}^*(s) \end{pmatrix} = \begin{pmatrix} 0 & 1 & 0 \\ -\cosh\theta & 0 & \sinh\theta \\ -\sinh\theta & 0 & \cosh\theta \end{pmatrix} \begin{pmatrix} \vec{t}(s) \\ \vec{n}(s) \\ \vec{b}(s) \end{pmatrix}. \tag{2.17}$$

2) If $\vec{D}$ is a timelike vector, then

$$\begin{pmatrix} \vec{t}^*(s) \\ \vec{n}^*(s) \\ \vec{b}^*(s) \end{pmatrix} = \begin{pmatrix} 0 & 1 & 0 \\ \sinh\theta & 0 & -\cosh\theta \\ -\cosh\theta & 0 & \sinh\theta \end{pmatrix} \begin{pmatrix} \vec{t}(s) \\ \vec{n}(s) \\ \vec{b}(s) \end{pmatrix}. \tag{2.18}$$

The trace of an $\vec{X}$ oriented line along space curve $\vec{\gamma}(s)$ is generally a trajectory ruled surface. A parametric equation of this trajectory ruled surface generated by $\vec{X}$ oriented line is given by

$$\varphi(s,k) = \vec{\gamma}(s) + k\vec{d}(s), \quad s,k \in I \subset \mathbb{R},$$

where $\vec{d}$ is the unit direction vector of $\vec{X}$ oriented line. The trace of an $\vec{X}$ oriented line along a space curve $\vec{\gamma}(s)$ is generally a trajectory ruled surface. The distribution parameter (or drall) of the $\varphi(s,k)$ trajectory ruled surface is given as

$$\delta_d = \frac{\det\left(\dot{\vec{\gamma}}, \vec{d}, \dot{\vec{d}}\right)}{\left\|\dot{\vec{d}}\right\|^2}. \tag{2.19}$$

A developable trajectory ruled surface is characterized by $\delta_d = 0$. If there exist a common perpendicular to two constructive rulings in the ruled surface, then the foot of the common perpendicular on the main ruling is called a central point. The locus of the central point is called striction curve. The parametrization of the sitriction curve on a trajectory ruled surface in [12] is given

$$\vec{C}(s)=\vec{\gamma}(s)-\frac{\left\langle \vec{\dot{\gamma}},\vec{\dot{d}}\right\rangle}{\left\|\vec{\dot{d}}\right\|^2}\vec{d}\,. \tag{2.20}$$

**3. Involute Trajectory Timelike Ruled Surfaces**

Let $\vec{\gamma}(s)$ be the spacelike involutes of a timelike space curve $\vec{r}(s)$ with spacelike $\vec{D}$ and $\left\{\vec{t}^*(s),\vec{n}^*(s),\vec{b}^*(s)\right\}$ be its Frenet frame defined as in Eq. (2.17). And also, we consider that a timelike oriented line $\vec{X}(s)$ in $\mathbb{R}_1^3$ such that it is firmly connected to Frenet frame of the involutes $\vec{\gamma}(s)$ is represented, uniquely with respect to this frame, in the form

$$\vec{X}(s)=x_1\vec{t}^*(s)+x_2\vec{n}^*(s)+x_3\vec{b}^*(s),\ \left\langle \vec{X},\vec{X}\right\rangle<0,\ \left\|\vec{X}\right\|=1. \tag{3.1}$$

where $x_i$ (i=1, 2, 3) are scalar of the arc length parameter of the involutes $\vec{\gamma}(s)$. The trajectory ruled surfaces generated by line $\vec{X}(s)$, $\vec{t}^*(s)$, $\vec{n}^*(s)$ and $\vec{b}^*(s)$ are

$$M:\varphi(s,v)=\vec{\gamma}(s)+v\vec{X}(s), \tag{3.2}$$

$$M_1:\varphi(s,u)=\vec{\gamma}(s)+u\vec{t}^*(s), \tag{3.3}$$

$$M_2:\varphi(s,z)=\vec{\gamma}(s)+z\vec{n}^*(s), \tag{3.4}$$

$$M_3:\varphi(s,w)=\vec{\gamma}(s)+w\vec{b}^*(s), \tag{3.5}$$

respectively. We can obtained the distribution parameter of the involute trajectory timelike ruled surface generated by $\vec{X}(s)$ in $\mathbb{R}_1^3$. Analytically, from equations (2.17), (3.1) and Frenet formulas

$$\begin{aligned}\vec{\dot{X}}&=\left(x_1\kappa-\dot{\theta}x_2\sinh\theta-\dot{\theta}x_3\cosh\theta\right)\vec{t}-x_2\left\|\vec{D}\right\|\vec{n}+\\&\quad\left(-x_1\tau+\dot{\theta}x_2\cosh\theta-\dot{\theta}x_3\sinh\theta\right)\vec{b}\end{aligned} \tag{3.6}$$

By differentiating Eq. (2.14) with respect to the arc length parameter *s*, we have

$$\vec{\dot{\gamma}}(s)=(c-s)\kappa\vec{n} \tag{3.7}$$

By substituting Eqs. (3.6) and (3.7) into Eq. (2.19), the distribution parameter of this surface is

$$\delta_{\vec{X}} = \frac{\det\left(\dot{\vec{\gamma}}, \vec{X}, \dot{\vec{X}}\right)}{\left\| \dot{\vec{X}} \right\|^2}, \tag{3.8}$$

$$\delta_{\vec{X}} = \frac{(c-s)\left[x_1x_3\left\|\vec{D}\right\| - \dot{\theta}\left(x_3^2 - x_2^2\right)\right]}{\left|\left(x_2^2 - x_1^2\right)\left\|\vec{D}\right\|^2 + \left(x_2^2 - x_3^2\right)\dot{\theta}^2 + 2x_1x_3\dot{\theta}\left\|\vec{D}\right\|\right|}, \tag{3.9}$$

The ruled surface developable if and only if $\delta_{\vec{X}}$ is zero. From Eq. (3.9) we have

$$(c-s)\left[x_1x_3\|\omega\| - \dot{\theta}\left(x_3^2 - x_2^2\right)\right] = 0 \tag{3.10}$$

Thus, we state the following theorem.

**Theorem 1.** The involute trajectory timelike ruled surface M is developable if and only if the Lorentzian timelike angle $\theta$ between $\vec{t}$ and $\vec{c}$ of space curve $\vec{r}(s)$ satisfies the following equality

$$\theta = \frac{x_1x_3}{x_3^2 - x_2^2}\int\left\|\vec{D}\right\|ds + \lambda, \tag{3.11}$$

where $\lambda$ is an arbitrary constant.

## 4. Special Cases

### 4.1 The Case $\vec{X}(s) = \vec{t}^*(s)$

In this case $x_1 = 1,\ x_2 = x_3 = 0$. Thus from Eq. (3.9)

$$\delta_{\vec{t}^*} = 0. \tag{4.1}$$

Thus we can give the following result.

**Result 1.** The ruled surface $M_1$ given by Eq. (3.3) is developable.

**4.2 The Case $\vec{X}(s)=\vec{n}^*(s)$**

In this case $x_2=1,\ x_1=x_3=0.$ Thus from Eq. (3.9)

$$\delta_{\vec{n}^*}=\frac{(c-s)\kappa\dot{\theta}}{\left|\left\|\vec{D}^2\right\|+\dot{\theta}^2\right|}. \tag{4.2}$$

So we can give the following result.

**Result 2.** If $\theta=$ constant (i.e. space curve $\vec{r}(s)$ is a general helix) or $c=s$ (i.e. space curve $\vec{r}(s)$ is coincident with $\vec{\gamma}(s)$) then the ruled surface $M_2$ given by Eq. (3.4) is developable.

**4.3 The Case $\vec{X}(s)=\vec{b}^*(s)$**

In this case $x_3=1,\ x_1=x_2=0$ and from Eq. (3.9)

$$\delta_{\vec{b}^*}=\frac{(c-s)\kappa}{\left|\dot{\theta}\right|}. \tag{4.3}$$

So we can give the following result.

**Result 3.** If $c=s$ (i.e. space curve $\vec{r}(s)$ is coincident with $\vec{\gamma}(s)$) then the ruled surface $M_3$ given by Eq. (3.5) is developable.

From the Eqs. (4.2) and (4.3) we have

$$\frac{\delta_{\vec{n}^*}}{\delta_{\vec{b}^*}}=\frac{\dot{\theta}^2}{\left\|\vec{D}^2\right\|+\dot{\theta}^2}. \tag{4.4}$$

Thus, the following theorem can be given.

**Theorem 2.** If $\vec{\gamma}(s)$ is the spacelike involutes of a timelike space curve $\vec{r}(s)$ with spacelike $\vec{D}$ and $\theta$ is the Lorentzian timelike angle between $\vec{t}$ and $\vec{c}$ then there is the relationship (4.4) between Darboux vector of $\vec{r}(s)$ and the distribution parameters of the ruled surfaces generated by $\vec{n}^*$ and $\vec{b}^*$ in $\mathbb{R}^3_1$.

From Eq. (4.4), if $\theta$ is constant then $\frac{\delta_{\vec{n}^*}}{\delta_{\vec{b}^*}}=0$. On the contrary, if $\frac{\delta_{\vec{n}^*}}{\delta_{\vec{b}^*}}=0$ then $\theta$ is constant. Therefore, with respect to this condition, we can give the following theorem.

**Theorem 3.** The space curve $\vec{r}(s)$ is general helix if and only if $\dfrac{\delta_{\vec{n}^*}}{\delta_{\vec{b}^*}} = 0.$

**4.4 The Case $\vec{X}(s)$ is in the Normal Plane**

In this case $x_1$ is zero $\left(-x_2^{\,2}+x_3^{\,2}=1\right)$. From Eq. (3.9), the distribution parameter of the involute trajectory timelike ruled surface $M$ is

$$\delta_{\vec{X}} = \frac{-(c-s)\kappa\dot{\theta}}{\left|x_2^{\,2}\left\|\vec{D}\right\|^2-\dot{\theta}^2\right|}. \tag{4.5}$$

Hence we state the following result.

**Result 4.** If the space curve $\vec{r}(s)$ is a general helix or coincident with $\vec{\gamma}(s)$ then the involute trajectory timelike ruled surface M generated by a timelike oriented line $\vec{X}(s)$ is developable on the normal plane in $\mathbb{R}^3_1$.

**4.5 The Case $\vec{X}(s)$ is in the Osculating Plane**

In this case $x_3$ is zero $\left(x_1^{\,2}-x_2^{\,2}=1\right)$. From Eq. (3.9), the distribution parameter of the involute trajectory timelike ruled surface $M$ is

$$\delta_{\vec{X}} = \frac{(c-s)\kappa\dot{\theta}x_2^{\,2}}{\left|-\left\|\vec{D}\right\|^2+x_2^{\,2}\dot{\theta}\right|}. \tag{4.6}$$

If $x_2$ is zero then $\vec{X}=\vec{t}^*$. This is the case 1.If $\theta=$ constant then $\vec{r}(s)$ is a general helix. If $c=s$ then $\vec{r}(s)$ is coincident with $\vec{\gamma}(s)$. Therefore Result 4. can be restated for the osculating plane in $\mathbb{R}^3_1$.

**4.6 The Case $\vec{X}(s)$ is in the Rectifying Plane**

In this case $x_2$ is zero $\left(x_1^{\,2}+x_3^{\,2}=1\right)$. From Eq. (3.9), the distribution parameter of the involute trajectory timelike ruled surface $M$ is

$$\delta_{\vec{X}} = \frac{(c-s)\left[x_1x_3\left\|\vec{D}\right\| - \dot{\theta}x_3^{\,2}\right]}{\left|-x_1^{\,2}\left\|\vec{D}\right\|^2 - x_3^{\,2}\dot{\theta}^2 + 2x_1x_3\dot{\theta}\left\|\vec{D}\right\|\right|}, \tag{4.7}$$

If $c = s$ then $\delta_{\vec{X}} = 0$. In this case, $\vec{r}(s)$ is coincident with $\vec{\gamma}(s)$. Also we can give the following theorem to characterize developability on the rectifying plane in $\mathbb{R}_1^3$.

**Theorem 4.** The involute trajectory timelike ruled surface M is developable on the rectifying plane if and only if the Lorentzian timelike angle $\theta$ between $\vec{t}$ and $\vec{c}$ of space curve $\vec{r}(s)$ satisfies the following equality

$$\theta = \frac{x_1}{x_3}\int\left\|\vec{D}\right\| ds + \lambda, \tag{4.8}$$

where $\lambda$ is an arbitrary constant.

### 4.7 The Case the base curve $\vec{\gamma}(s)$ is the striction curve $\vec{C}(s)$

From Eq. (2.20) the parametrization of the striction curve on the involute trajectory timelike ruled surface generated by timelike oriented line $\vec{X}(s)$ is given by

$$\vec{C}(s) = \vec{\gamma}(s) + \frac{x_2(c-s)\kappa\left\|D\right\|}{\left\|\dot{\vec{X}}\right\|^2}\vec{X}. \tag{4.9}$$

From here if the base curve $\vec{\gamma}(s)$ is the striction curve $\vec{C}(s)$ then we have $x_2 = 0$ or $c = s$. Hence the following theorem can be given.

**Result 5.** If the base curve $\vec{\gamma}(s)$ is the same as the striction curve $\vec{C}(s)$ then the oriented line $\vec{X}(s)$ of the involute trajectory timelike ruled surface is on the rectifying plane or space curve $\vec{r}(s)$ is coincident with $\vec{\gamma}(s)$.

**Example:** Let $r(s)=\left(2\sinh\left(\frac{s}{\sqrt{3}}\right),2\cosh\left(\frac{s}{\sqrt{3}}\right),\frac{s}{\sqrt{3}}\right)$ be a unit speed timelike helix with spacelike Darboux vector $\vec{D}$ such that $\kappa=\frac{2}{3}$ and $\tau=\frac{1}{3}$. The short calculations give

$$\begin{cases}\vec{t}(s)=\left(\frac{2}{\sqrt{3}}\cosh\left(\frac{s}{\sqrt{3}}\right),\frac{2}{\sqrt{3}}\sinh\left(\frac{s}{\sqrt{3}}\right),\frac{1}{\sqrt{3}}\right)\\ \vec{n}(s)=\left(\sinh\left(\frac{s}{\sqrt{3}}\right),\cosh\left(\frac{s}{\sqrt{3}}\right),0\right)\\ \vec{b}(s)=\left(\frac{1}{\sqrt{3}}\cosh\left(\frac{s}{\sqrt{3}}\right),\frac{1}{\sqrt{3}}\sinh\left(\frac{s}{\sqrt{3}}\right),\frac{2}{\sqrt{3}}\right)\end{cases}. \tag{4.10}$$

In this situation, from Eq. (2.19) the involutes $\vec{\gamma}(s)$ of the curve $\vec{r}(s)$ can be given by the equation

$$\begin{aligned}\gamma(s)=\Bigg(&2\sinh\left(\frac{s}{\sqrt{3}}\right)+|c-s|\frac{2}{\sqrt{3}}\cosh\left(\frac{s}{\sqrt{3}}\right),\\ &2\cosh\left(\frac{s}{\sqrt{3}}\right)+|c-s|\frac{2}{\sqrt{3}}\sinh\left(\frac{s}{\sqrt{3}}\right),\frac{c}{\sqrt{3}}\Bigg).\end{aligned} \tag{4.11}$$

From Eqs. (2.9) and (2.10) we have

$$\vec{D}=\frac{1}{3}\vec{t}-\frac{2}{3}\vec{b} \tag{4.12}$$

$$\begin{cases}\cosh\theta=\frac{2\sqrt{3}}{3}\\ \sinh\theta=\frac{\sqrt{3}}{3}\end{cases} \tag{4.13}$$

respectively. By using Eq. (2.17) we have the Frenet trihedron of the involutes $\vec{\gamma}(s)$ of the curve $\vec{r}(s)$

$$\begin{cases} \vec{t}^{\,*}(s) = \left( \sinh\left(\dfrac{s}{\sqrt{3}}\right), \cosh\left(\dfrac{s}{\sqrt{3}}\right), 0 \right) \\ \vec{n}^{\,*}(s) = \left( -\cosh\left(\dfrac{s}{\sqrt{3}}\right), -\sinh\left(\dfrac{s}{\sqrt{3}}\right), 0 \right). \\ \vec{b}^{\,*}(s) = (0,0,1) \end{cases} \tag{4.14}$$

Thus we obtain the involute trajectory ruled surfaces generated by $\vec{t}^{\,*}$, $\vec{n}^{\,*}$ and $\vec{b}^{\,*}$ as

$$\begin{aligned} \varphi_{\vec{t}^{\,*}}(s,v) = \Bigg( & 2\sinh\left(\frac{s}{\sqrt{3}}\right) + |c-s|\frac{2}{\sqrt{3}}\cosh\left(\frac{s}{\sqrt{3}}\right) + v\sinh\left(\frac{s}{\sqrt{3}}\right), \\ & 2\cosh\left(\frac{s}{\sqrt{3}}\right) + |c-s|\frac{2}{\sqrt{3}}\sinh\left(\frac{s}{\sqrt{3}}\right) + v\cosh\left(\frac{s}{\sqrt{3}}\right), \frac{c}{\sqrt{3}} \Bigg), \end{aligned} \tag{4.15}$$

$$\begin{aligned} \varphi_{\vec{n}^{\,*}}(s,v) = \Bigg( & 2\sinh\left(\frac{s}{\sqrt{3}}\right) + |c-s|\frac{2}{\sqrt{3}}\cosh\left(\frac{s}{\sqrt{3}}\right) - v\cosh\left(\frac{s}{\sqrt{3}}\right), \\ & 2\cosh\left(\frac{s}{\sqrt{3}}\right) + |c-s|\frac{2}{\sqrt{3}}\sinh\left(\frac{s}{\sqrt{3}}\right) - v\sinh\left(\frac{s}{\sqrt{3}}\right), \frac{c}{\sqrt{3}} \Bigg), \end{aligned} \tag{4.16}$$

$$\begin{aligned} \varphi_{\vec{b}^{\,*}}(s,v) = \Bigg( & 2\sinh\left(\frac{s}{\sqrt{3}}\right) + |c-s|\frac{2}{\sqrt{3}}\cosh\left(\frac{s}{\sqrt{3}}\right), \\ & 2\cosh\left(\frac{s}{\sqrt{3}}\right) + |c-s|\frac{2}{\sqrt{3}}\sinh\left(\frac{s}{\sqrt{3}}\right), \frac{c}{\sqrt{3}} + v \Bigg), \end{aligned} \tag{4.17}$$

respectively, where $0 \le s \le \pi$, $-2 \le v \le 2$ and $c = 1$ (Figs. 1-3).

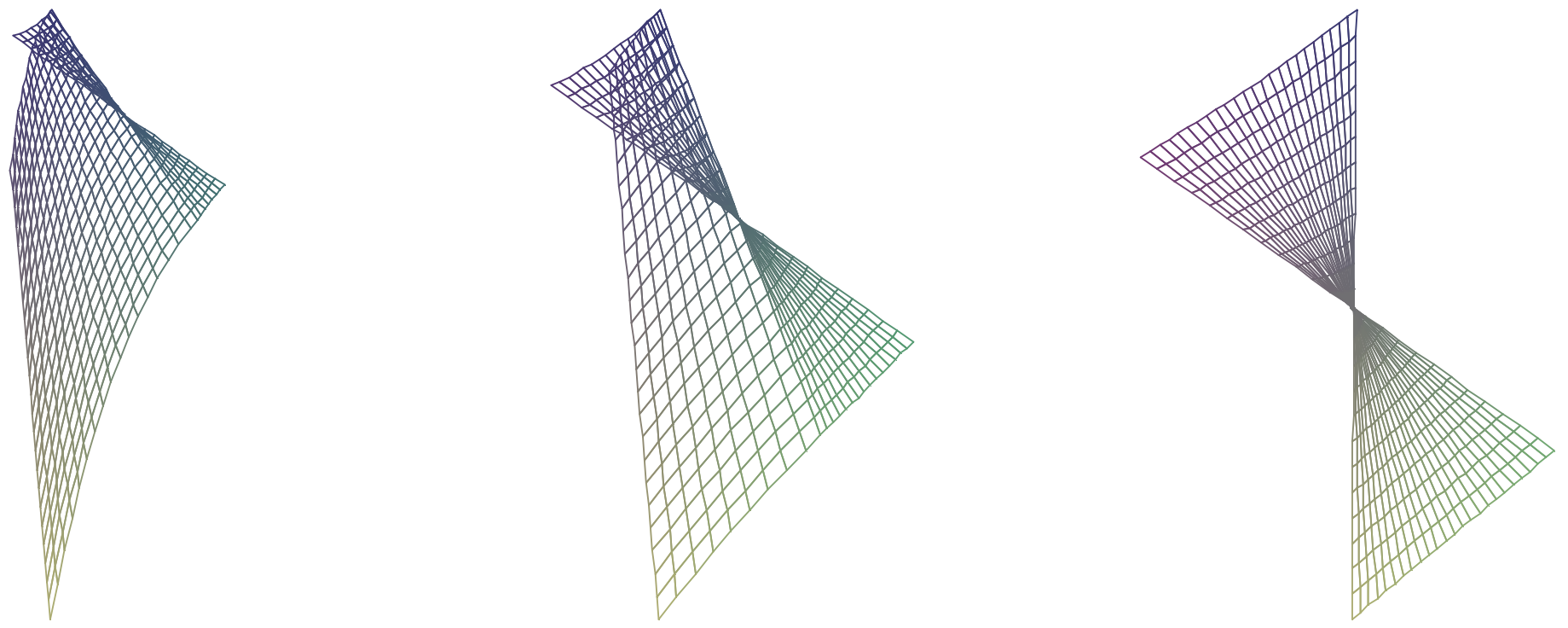

**Fig. 1** Some involute trajectory ruled surfaces generated by $\vec{t}^*$

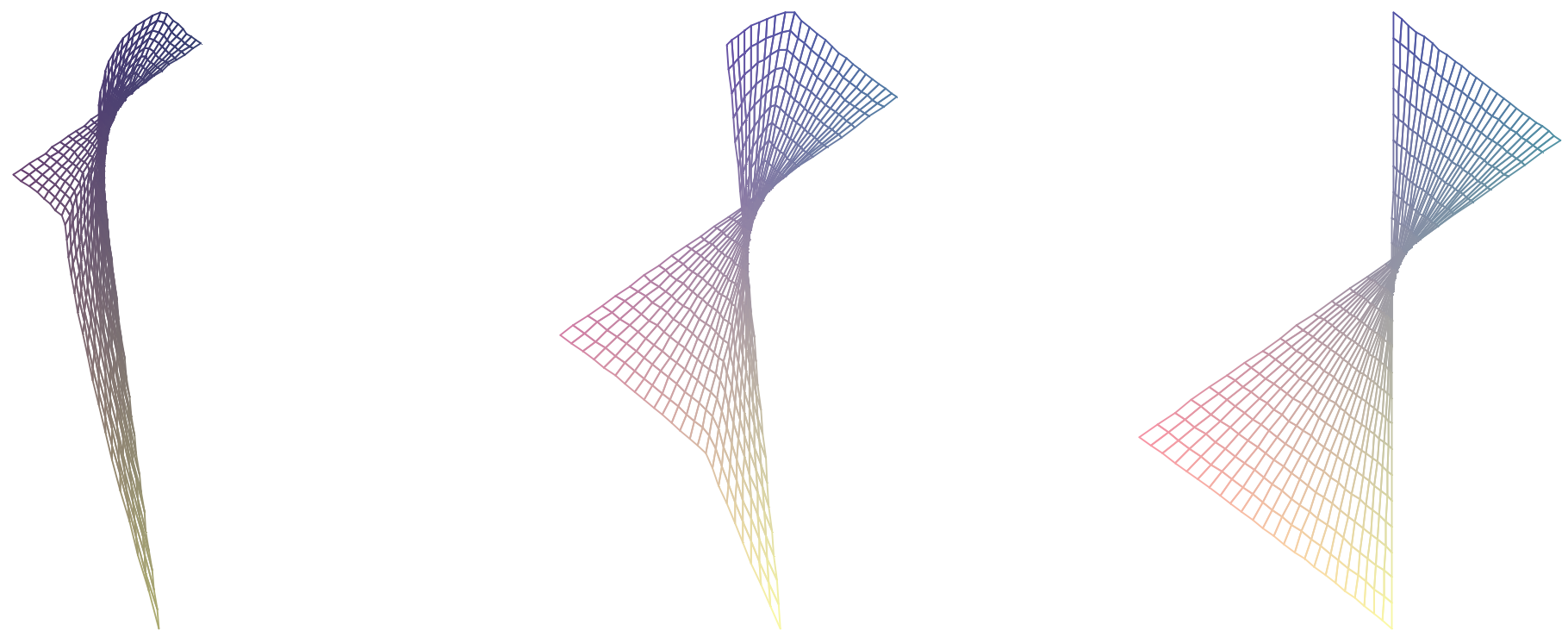

**Fig. 2** Some involute trajectory ruled surfaces generated by $\vec{n}^*$

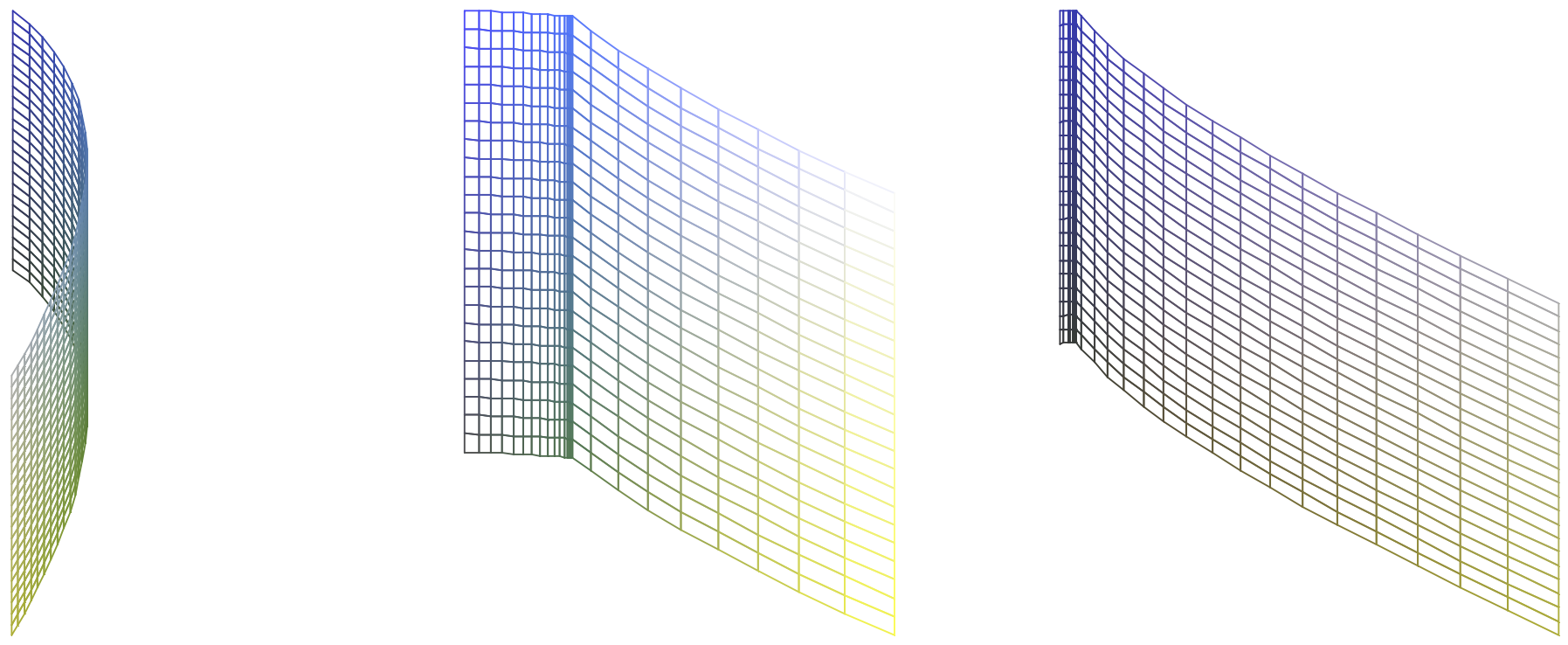

**Fig. 3** Some involute trajectory ruled surfaces generated by $\vec{b}^*$

## Appendix

The same procedure is also applied for the spacelike involutes $\vec{\gamma}(s)$ of a timelike space curve $\vec{r}(s)$ with timelike $\vec{D}$ given by Eq. (2.18), it is seen that the same theorems and results can be easily obtained.